\documentclass[11pt,a4paper]{article}
\usepackage{bbm}
\usepackage{mathrsfs}
\usepackage{amsfonts}
\usepackage{}

\usepackage{graphicx}
\usepackage{amsfonts,amsmath, amssymb}
\usepackage{color}

\newtheorem{theo}{\textbf{\ \ \quad Theorem}}[section]

\newtheorem{lem}{\textbf{\ \ \quad Lemma}}[section]
\newtheorem{remark}{\textbf{\ \ \quad Remark}}[section]
\newtheorem{col}{\textbf{\ \ \quad Corollary}}[section]

\newcommand{\lbl}[1]{\label{#1}}

\newcommand{\be}{\begin{equation}}
\newcommand{\ee}{\end{equation}}
\newcommand\bes{\begin{eqnarray}}
\newcommand\ees{\end{eqnarray}}
\newcommand{\bess}{\begin{eqnarray*}}
\newcommand{\eess}{\end{eqnarray*}}

\newcommand{\nm}{\nonumber}

{\setlength\arraycolsep{2pt}}

\title{Regularity of solutions to elliptic equations with Grushin's operator}

\author{Xiaohuan Wang$^{a,b}$ and Jihui Zhang$^a$\\
\\
\ \\
    {\small \it $^a$ Institute of Mathematics, School of Mathematical Science}\\
  {\small \it Nanjing Normal University, Nanjing 210023, China}\\
  {\small \it $^a$ School of Mathematics and Statistics, Henan University}\\
  {\small \it Kaifeng, Henan 475001, China}\\
  {\small \tt xiaohuanw@126.com}
}

\begin{document}
\maketitle

\medskip

\begin{abstract} In this paper, we consider the regularity of solutions to
elliptic equation with Grushin's operator. By using the Feynman-Kac formula,
we first get the expression of heat kernel, and then by using the properties of
heat kernel, the optimal regularity of solutions will be obtained.
The novelty of this paper is that the Grushin's operator is a
degenerate operator.



\end{abstract}

\baselineskip=15pt

\section{Introduction}
\setcounter{equation}{0}

The regularity of solutions to second order elliptic equation has been extensively
studied by many authors, see the book \cite{GT1983}. But for degenerate
elliptic equation, there is few work about  the regularity results  until now.
In this paper, we focus on a special degenerate elliptic equation--Grushin's elliptic equation. The
main reason why we can deal with it is that we can get the expression of heat kernel
by using the probability method.

For Grushin's operator, many authors studied it.
 Beckner \cite{Be2001} obtained the Sobolev estimates for the
 Grushin's operator in low dimensions by using
 hyperbolic symmetry and conformal geometry. Riesz transforms and
multipliers for the Grushin's operator was considered by Jotsaroop et al.
\cite{JST2013}. Tri \cite{Tri1998} studied the generalized  Grushin's equation.
$L^p$-estimates for the wave equation
associated to the Grushin operator was studied by Jotsaroop-Thangavelu
\cite{JT2012}. The fundamental solution for a degenerate parabolic pseudo-differential operator covering
the Grushin's operator was obtained by Tsutsumi \cite{S1974,S1975}. Furthermore,
Tsutsumi \cite{S1976} constructed a left parametrix for a pesudo-differential
operator. We remark that Tsutsumi did not give the exact expression of heat kernel.
In the book \cite{CCFI2009}, they gave the expression of heat kernel (see Page 191), but
he expression is hard to use because they used the inverse Fourier transform.
Similar degenerate elliptic equation was studied by
Robinson-Sikora \cite{RS2008} and the Hardy inequalities for
Grushin's operator was considered by \cite{YSK2015}.

In this paper, in view of probability point, we give a new expression and then get
the regularity of solution by using the properties of heat kernel.

This paper is arranged as follows. In next section, some
preliminaries are given and the main results will be proved in section 3.
Throughout this paper, we write $C$ as a general positive constant and $C_i$, $i=1,2,\cdots$ as
a concrete positive constant.

\section{Main results}
\setcounter{equation}{0}

Consider the Grushin's operator
   \bess
\mathcal{L}=\frac{1}{2}(\partial_{x}^2+x^2\partial_{y}^2),
   \eess
which is the generator of the diffusion process $(X_t,Y_t)$, where $(X_t,Y_t)$ satisfies
  \bess\left\{\begin{array}{llll}
dX_t=dW^1_t,\\
dY_t=X_tdW^2_t,\\
X_0=\mu_1,\ \ \ Y_0=\mu_2.
  \end{array}\right.\eess
Here $W^i_t$ $(i=1,2)$ are standard i.i.d Brownian motion. It is easy to see that
the process $(X_t,Y_t)$ is a Gaussian stochastic process. Direct calculations
show that
   \bess
\mathbb{E}\left(\begin{array}{cccc}  X_t\\[-.5mm]  Y_t
\end{array}\right)=\left(\begin{array}{cccc}  \mu_1\\[-.5mm]  \mu_2
\end{array}\right),\ \ \
Cov(X_t,Y_t)=\left(\begin{array}{cccc}
t &\ \mu_1t & \\[-.5mm]
\mu_1t &\ \mu_1^2t+\frac{1}{2}t^2 &
\end{array}\right) .
   \eess
Therefore, we get the heat kernel of the
operator $\mathcal{L}$ is
  \bess
K(t,x,\mu_1,y,\mu_2)=\frac{1}{2\pi t^{3/2}}\exp\left\{-\frac{(x-\mu_1)^2}{t}-\frac{[\mu_1(x-\mu_1)-y+\mu_2]^2}{t^2}\right\},
   \eess
which yields that
  \bess
\nabla_xK(t,x,y)&=&-\frac{2x}{t}K(t,x,y),\\
\nabla_yK(t,x,y)&=&-\frac{y}{t^2}K(t,x,y).
  \eess

 We first
want to solve the following elliptic equation
   \bes
\mathcal{L^\lambda}_b:=(\mathcal{L}-\lambda)u+b\cdot\nabla u=f.
   \lbl{0.1}\ees
   \begin{theo}\lbl{t1.1} Assume that $b\in C^\beta_b(\mathbb{R}^2)$. There exists
a $\lambda_0>0$ such that for any $f\in C^\beta_b(\mathbb{R}^2)$ and $\lambda>\lambda_0$, there is
a unique solution $u\in C_b^{1+\beta}(\mathbb{R}^2)$ to equation (\ref{0.1}) such that
   \bes
\lambda^{\delta}\|u\|_{C_b^{1+\beta}(\mathbb{R}^2)}\leq\|f\|_{ C_b^{\beta}(\mathbb{R}^2)},
   \lbl{0.2}\ees
where $\delta>0$ is defined in Lemma \ref{l2.1}.

Assume further that $\nabla_yh(x,\cdot)\in C_b^\beta(\mathbb{R})$ for any $x\in\mathbb{R}$, where $h=b$ or $f$.
Then there is
a unique solution $u\in C_b^{1+\beta}(\mathbb{R}^2)$ to equation (\ref{0.1}) such that
   \bess
\lambda^{\delta}\|u\|_{C_b^{2+\beta}(\mathbb{R}^2)}\leq\|f\|_{ C_b^{\beta}(\mathbb{R};C^{1+\beta}(\mathbb{R}))}.
    \eess
  \end{theo}

Let
   \bess
u(x,y):=(\lambda-\mathcal{L})^{-1}f(x)=\int_0^\infty e^{-\lambda t}K(t,\cdot)\ast f(x,y)dt.
  \eess
Next, we consider the $L^p$-regularity of $u(x,y)$.
\begin{theo}\lbl{t1.2}
 Assume that $f(\cdot,y)\in L^{p}(\mathbb{R})$ for every $y\in\mathbb{R}$ with  and
$f(x,\cdot)\in W^{s,q}(\mathbb{R},\mathbb{R})$ for every $y\in\mathbb{R}$, $0<s<1$ and
$q>1$. Then we have the following estimates:
  \bess
\|u\|_{L^r(\mathbb{R}^2)}&\leq&C\lambda^{-1-\frac{3}{2r}+\frac{1}{2p}+\frac{1}{q}}\|f\|_{L^p(\mathbb{R},L^q(\mathbb{R}))},\\
\|\nabla_xu\|_{L^r(\mathbb{R}^2)}&\leq&C\lambda^{-\frac{1}{2}-\frac{3}{2r}+\frac{1}{2p}+\frac{1}{q}}\|f\|_{L^p(\mathbb{R},L^q(\mathbb{R}))},\\
\|\nabla_yu\|_{L^r(\mathbb{R}^2)}&\leq&C\lambda^{-s-\frac{3}{2r}+\frac{1}{2p}}\|f\|_{L^p(\mathbb{R},W^{s,q}(\mathbb{R}))}.\\
  \eess
Moreover, if we take $0<s<1$ and $p,q,r>1$ such that $-\frac{3}{2r}+\frac{1}{2p}+\frac{1}{q}\leq0$
and $-s-\frac{3}{2r}+\frac{1}{2p}+1\leq0$, we have
 \bess
 \|\nabla^2_xu\|_{L^r(\mathbb{R}^2)}\leq C,\ \ \
\|\nabla^2_yu\|_{L^r(\mathbb{R}^2)}\leq C.
   \eess
\end{theo}

 \section{Proof of Main results}
\setcounter{equation}{0}

Denote $[\cdot]_\beta$ by the semi-norm of $C^\beta$.
Set
   \bess
\mathcal{A}:=\left\{f(x,y):\begin{array}{lllll} f(\cdot,y)\in L^\infty(\mathbb{R})\  for\  any \ y\in\mathbb{R},
f(x,\cdot)\in C_b^\beta(\mathbb{R})\\   for \ any \
 x\in\mathbb{R}\ with \ 0<\beta<1\end{array}\right\}
   \eess

\begin{lem}\lbl{l2.1} Assume that $f\in \mathcal{A}$. Then
  \bess
\|u\|_{L^\infty(\mathbb{R}^2)} &\leq& C\lambda^{-\frac{1}{2}}\|f\|_{L^\infty(\mathbb{R}^2)},\ \qquad \ \
\|\nabla_xu\|_{L^\infty(\mathbb{R}^2)}\leq C\lambda^{-\frac{1}{2}}\|f\|_{L^\infty(\mathbb{R}^2)},\\
\|\nabla_yu\|_{L^\infty(\mathbb{R}^2)}&\leq&C\lambda^{-\beta}\|f\|_{L^\infty(\mathbb{R},C_b^\beta(\mathbb{R}))},\ \  \ \
[\nabla_xu]_\beta\leq C\lambda^{-\frac{1}{2}+\beta}\|f\|_{L^\infty(\mathbb{R},C_b^\beta(\mathbb{R}))}.
  \eess
Moreover, if $f(\cdot,y)\in C_b^\beta(\mathbb{R})$ for any
$y\in \mathbb{R}$, it holds that
  \bess
[\nabla_yu]_\beta\leq C\lambda^{-\delta}\|f\|_{C^\beta(\mathbb{R}^2)},
  \eess
where $0<\delta<(\beta\wedge(1-\beta))$. That is to say, $\nabla u\in C_b^\beta(\mathbb{R}^2)$.
  \end{lem}

{\bf Proof.}
Simply calculations show that
   \bess
\|u\|_{L^\infty(\mathbb{R}^2)}&=&\|\int_0^\infty e^{-\lambda t}K(t,\cdot)\ast f(x,y)dt\|_{L^\infty(\mathbb{R}^2)}\\
&\leq&\int_0^\infty e^{-\lambda t}\|K(t,\cdot)\|_{L^1(\mathbb{R})}\|f\|_{L^\infty(\mathbb{R}^2)}dt\\
&\leq&C\lambda^{-1}\|f\|_{L^\infty(\mathbb{R}^2)},
  \eess
   \bess
\|\nabla_xu\|_{L^\infty(\mathbb{R}^2)}&=&\|\int_0^\infty e^{-\lambda t}\nabla_xK(t,\cdot)\ast f(x,y)dt\|_{L^\infty(\mathbb{R}^2)}\\
&\leq&\int_0^\infty e^{-\lambda t}\|\nabla_xK(t,\cdot)\|_{L^1(\mathbb{R})}\|f\|_{L^\infty(\mathbb{R}^2)}dt\\
&\leq&C\lambda^{-\frac{1}{2}}\|f\|_{L^\infty(\mathbb{R}^2)},
  \eess
and
    \bess
&&\|\nabla_yu\|_{L^\infty(\mathbb{R}^2)}\\
&=&\|\int_0^\infty e^{-\lambda t}\nabla_yK(t,\cdot)\ast f(x,y)dt\|_{L^\infty(\mathbb{R}^2)}\\
&=&\|\int_0^\infty e^{-\lambda t}\int_{\mathbb{R}^2}\nabla_yK(t,x-u,y-v)(f(u,v)-f(u,y))dudvdt\|_{L^\infty(\mathbb{R}^2)}\\
&\leq&\|f\|_{L^\infty(\mathbb{R},C_b^\beta(\mathbb{R}))}\|\int_0^\infty e^{-\lambda t}\int_{\mathbb{R}^2}|\nabla_yK(t,x-u,y-v)|\cdot|y-v|^\beta dudv dt\|_{L^\infty(\mathbb{R}^2)}\\
&\leq&C\|f\|_{L^\infty(\mathbb{R},C_b^\beta(\mathbb{R}))}\int_0^\infty e^{-\lambda t}t^{-1+\beta}dt\\
&\leq&C\lambda^{-\beta}\|f\|_{L^\infty(\mathbb{R},C_b^\beta(\mathbb{R}))},
  \eess
where we used the fact that
   \bess
\int_{\mathbb{R}^2}\nabla_yK(t,x-u,y-v)f(u,y)dudv=\int_{\mathbb{R}}\left(\int_{\mathbb{R}}\nabla_yK(t,x-u,y-v)dv\right)f(u,y)du=0.
   \eess
We also remark that the meaning of $\|f\|_{L^\infty(\mathbb{R},C_b^\beta(\mathbb{R}))}$ is that we take
infinity norm for the first variable and take H\"{o}lder norm for the second variable.

Recall the following interpolation inequality
   \bess
[u]_{\sigma\alpha+(1-\sigma)\gamma}\leq ([u]_{\alpha})^\sigma([u]_\gamma)^{1-\sigma}, \ \
0\leq\alpha<\gamma\leq1,\ \ \sigma\in(0,1).
  \eess
Now, if $0<\beta<1$, applying the above inequality with $\alpha=0,\gamma=1$ and
$\beta=\gamma(1-\sigma)$, we have
   \bess
[\nabla_xu]_\beta&=&\left[\int_0^\infty e^{-\lambda t}\nabla_xK(t,\cdot)\ast f(x,y)dt\right]_\beta\\
&\leq&\int_0^\infty e^{-\lambda t}\int_{\mathbb{R}^2}\left[\nabla_xK(t,\cdot,\cdot)\right]_\beta(x,y)dxdy\|f\|_{L^\infty(\mathbb{R}^2)}dt\\
&\leq&\|f\|_{L^\infty(\mathbb{R}^2)}\int_0^\infty t^{-\frac{1}{2}-\beta}e^{-\lambda t}dt\\
&\leq&C\lambda^{-\frac{1}{2}+\beta}\|f\|_{L^\infty(\mathbb{R}^2)}.
  \eess
Next, we consider the derivative of the second variable.
 \bess
[\nabla_yu]_\beta
&=&
\left[\int_0^\infty e^{-\lambda t}\nabla_yK(t,\cdot)\ast f(x,y)dt\right]_\beta\\
&\leq&\int_0^\infty e^{-\lambda t}\left(\sup_{x,\hat x,y,\hat y\in\mathbb{R},x\neq \hat x,y\neq\hat y}\frac{1}{|x-\hat x|^\beta+|y-\hat y|^\beta}\right.\\
&&\times\Big|\int_{\mathbb{R}^2}\nabla_yK(t,x-u,v)(f( u,y-v)-f(u,\hat y-v))dudv\\
&&\left.+
\int_{\mathbb{R}^2}\nabla_yK(t,u,\hat y-v)(f( x-u,v)-f(\hat x-u,v))dudv\Big|\right)dt\\
&=&\int_0^\infty e^{-\lambda t}\left(\sup_{x,\hat x,y,\hat y\in\mathbb{R},x\neq \hat x,y\neq\hat y}\frac{1}{|x-\hat x|^\beta+|y-\hat y|^\beta}|I_1+I_2|\right)dt.
  \eess
By dividing the real line into two parts, we have
   \bess
I_1&=&\int_{\mathbb{R}^2}\nabla_yK(t,x-u,y-v)(f( u,v)-f(u,y)dudv\\
&&+\int_{\mathbb{R}^2}\nabla_yK(t,x-u,\hat y-v)(f(u,\hat y)-f( u,v))dudv\\
&=&\int_{\mathbb{R}^2}\nabla_yK(t,u,y-v)(f(x-u,v)-f(x-u,y ))dudv\\
&&+\int_{\mathbb{R}^2}\nabla_yK(t,u,\hat y-v)(f(x-u,\hat y)-f(x-u,v))dudv\\
&=&\int_{\mathbb{R}}\int_{|y-v|\leq2|y-\hat y|}\nabla_yK(t,u, y-v)(f(x-u,v)-f(x-u, y)dudv\\
&&+\int_{\mathbb{R}}\int_{|y-v|\leq2|y-\hat y|}\nabla_yK(t,u,\hat y-v)(f(x-u,\hat y)-f(x-u,v))dudv\\
&&+\int_{\mathbb{R}}\int_{|y-v|>2|y-\hat y|}(\nabla_yK(t,u,y-v)-\nabla_yK(t,u,\hat y-v))(f(x-u,v)-f(x-u,y))dudv\\
&&+\int_{\mathbb{R}}\int_{|y-v|>2|y-\hat y|}\nabla_yK(t,u,\hat y-v)(f(x-u,\hat y)-f(x-u,y))dudv\\
&=:&I_{11}+\cdots+I_{14}.
  \eess
Let us estimate $I_{11}$-$I_{14}$. By using the form of heat kernel, we get
   \bess
|I_{11}|&=&\Big|\int_{\mathbb{R}}\int_{|y-v|\leq2|y-\hat y|}\nabla_yK(t,u, y-v)(f(x-u,v)-f(x-u,y-u^2)dudv\Big|\\
&\leq&Ct^{-\frac{3}{2}}\int_{\mathbb{R}}e^{-\frac{u^2}{t}}\left(\int_{|y-v|\leq2|y-\hat y|}\frac{|y-v|}{t^2}
e^{-\frac{(y-v)^2}{t^2}}|y-v|^\beta dv\right)du\\
&=&Ct^{-\frac{3}{2}}\int_{\mathbb{R}}e^{-\frac{u^2}{t}}\left(\int_{|z|\leq2|y-\hat y|}\frac{|z|}{t^2}
e^{-\frac{z^2}{t^2}}|z|^\beta dv\right)du\\
&\leq&Ct^{-1+\beta}|y-\hat y|^\beta.
   \eess
Similarly, we can obtain
    \bess
I_{12}\leq Ct^{-1+\beta}|y-\hat y|^\beta.
   \eess
Next, we consider $I_{13}$.
    \bess
|I_{13}|&=&\Big|\int_{\mathbb{R}}\int_{|y-v|>2|y-\hat y|}(\nabla_yK(t,u,y-v)-\nabla_yK(t,u,\hat y-v))(f(x- u,v)-f(x-u,y))dudv\Big|\\
&\leq&C\int_{\mathbb{R}}\int_{|y-v|>2|y-\hat y|}|\nabla_yK(t,u,y-v)-\nabla_yK(t,u,\hat y-v)|y-v|^\beta dudv\\
&\leq&Ct^{-\frac{3}{2}}\int_{\mathbb{R}}e^{-\frac{u^2}{t}}
\int_{|y-v|>2|y-\hat y|}\Big|\frac{y-v}{t^2}e^{-\frac{(y-v)^2}{t^2}}-
\frac{\hat y-v}{t^2}e^{-\frac{(\hat y-v)^2}{t^2}}\Big||y-v|^\beta dudv.
  \eess
Notice that $|y-v|>2|y-\hat y|$. So for every $\xi\in [y,\hat y]$,
  \bess
\frac{1}{2}|y-v| \leq |\xi-v|\leq 2|y-v|.
  \eess
We recall the following fractional mean value formula
(see (4.4) of \cite{Jumarie2006})
   \bess
f(x+h)=f(x)+\Gamma^{-1}(1+\beta)h^\gamma f^{(\gamma)}(x+\theta h),
   \eess
where $0<\gamma<1$ and $\theta>0$ depends on $h$ satisfying
   \bess
\lim\limits_{h\downarrow0}\theta^\gamma=\frac{\Gamma^2(1+\gamma)}{\Gamma(1+2\gamma)}.
   \eess
Denote
   \bess
\tilde K(t,v)=\frac{v}{t^2}e^{-\frac{v^2}{t^2}}.
   \eess
By using the above fractional mean value formula with $\gamma>\beta$
and the interpolation inequality in H\"{o}lder space
   \bess
|I_{13}|&\leq&Ct^{-\frac{3}{2}}|y-\hat y|^\beta\int_{\mathbb{R}}e^{-\frac{u^2}{t}}
\int_{|y-v|>2|y-\hat y|}[\tilde K(t,y-v)]_\gamma|y-\hat y|^{\gamma-\beta}|y-v|^\beta dudv\\
 &\leq&Ct^{-\frac{3}{2}}|y-\hat y|^\beta\int_{\mathbb{R}}e^{-\frac{u^2}{t}}du
\int_{|y-v|>2|y-\hat y|}[\tilde K(t,y-v)]_\gamma|y-v|^{\gamma} dv\\
&\leq&Ct^{-1+\gamma-\beta}|y-\hat y|^\beta\int_{\mathbb{R}}e^{-u^2}du\int_0^\infty|v|^\gamma e^{-\frac{1}{2}v^2}dv\\
&\leq&Ct^{-1+\gamma-\beta}|y-\hat y|^\beta.
  \eess
Lastly, by using the properties of heat kernel $K$, it is easy to see that
   \bess
&&\int_{|u^2-(y-v)|>2|y-\hat y|}\nabla_yK(t,u,y-v)dv\\
&=&\int_{|v|>2|y-\hat y|}\nabla_yK(t,u,u^2-v)dv\\
&=&(2\pi)^{-1}t^{-\frac{3}{2}}
e^{\frac{u^2}{t}}e^{-\frac{v^2}{t^2}}\Big|_{v=2|y-\hat y|}^{v=-2|y-\hat y|}\\
&=0&.
  \eess
Using the above equality and similar to the operation of $I_{13}$, we have
   \bess
|I_{14}|&=&\Big|\int_{\mathbb{R}}\int_{|y-v|>2|y-\hat y|}\nabla_yK(t,u,\hat y-v)(f( u,\hat y-u^2)-f(u,y-u^2))dudv\Big|\\
&=&\Big|\int_{\mathbb{R}}\int_{|y-v|>2|y-\hat y|}(\nabla_yK(t,u,\hat y-v)-\nabla_yK(t,u,y-v))\\
&&\times(f( u,\hat y-u^2)-f(u,y-u^2))dudv\Big|\\
&=&\Big|\int_{\mathbb{R}}(f(x- u,\hat y)-f(x-u,y))du\int_{|y-v|>2|y-\hat y|}\\
&&\times(\nabla_yK(t,u,\hat y-v)-\nabla_yK(t,u,y-v))dv\Big|\\
&\leq&C\|f\|_{L^\infty(\mathbb{R},C_b^\beta(\mathbb{R}))}|y-\hat y|^\beta\int_{\mathbb{R}}du\\
&&\times\int_{|y-v|>2|y-\hat y|}
|\nabla_yK(t,u,\hat y-v)-\nabla_yK(t,u,y-v)|dv\\
&\leq&
C\|f\|_{L^\infty(\mathbb{R},C_b^\beta(\mathbb{R}))}t^{-1+\gamma-\beta}|y-\hat y|^\beta.
  \eess
Substituting $I_{11}-I_{14}$ into $I_1$, we get
  \bess
&&\int_0^\infty e^{-\lambda t}\left(\sup_{x,\hat x,y,\hat y\in\mathbb{R},x\neq \hat x,y\neq\hat y}
\frac{1}{|x-\hat x|^\beta+|y-\hat y|^\beta}|I_1|\right)dt\\
&\leq&C\int_0^\infty e^{-\lambda t}\left(t^{-1+\beta}+t^{-1+\gamma-\beta}\right)dt\\
&\leq& C\lambda^{-(\beta\wedge(\gamma-\beta))}.
  \eess
Similarly, we can prove that if $f(\cdot,y)\in C_b^\beta(\mathbb{R})$ for any $y\in \mathbb{R}$,
  \bess
\int_0^\infty e^{-\lambda t}\left(\sup_{x,\hat x,y,\hat y\in\mathbb{R},x\neq \hat x,y\neq\hat y}
\frac{1}{|x-\hat x|^\beta+|y-\hat y|^\beta}|I_2|\right)dt\leq C\lambda^{-(\beta\wedge(\gamma-\beta))}.
  \eess
Summing the above discussion, we obtain
  \bess
[\nabla_yu]_\beta\leq C\lambda^{-(\beta\wedge(\gamma-\beta))}\|f\|_{C^\beta(\mathbb{R}^2)}.
  \eess
Noting that the above inequality holds for $0<\gamma<1$. The proof of this lemma is complete. $\Box$

\begin{remark}\lbl{r2.1} It is well known that if $f\in C^\alpha(\mathbb{R}^n)$, then the solution $u$ of the following
equation
   \bess
u_t-\Delta u=f,\ \ u_0=0,
   \eess
belongs to $C^{2+\alpha}(\mathbb{R}^n)$, which is the Schauder theory. Noting that the heat kernel of above
equation is Gauss heat kernel, that is
   \bess
K(t,x)=\frac{1}{(2\pi t)^{\frac{n}{2}}}e^{-\frac{x^2}{t}}.
   \eess
It is easy to see that $x\sim\sqrt{t}$. But in our case, different axis has different
scaling, that is,
   \bess
x\sim\sqrt{t},\ \ \ y\sim t.
  \eess
Thus when we take derivative for variable $x$, we can get $t^{-\frac{1}{2}}$, and
take double derivative for variable $x$, we will get $t^{-1}$. But if
we take derivative for variable $y$, we shall get $t^{-1}$, which is different from
the classical case. In Schauder theory, we can get $C^{2+\alpha}(\mathbb{R}^n)$ estimates,
but in our case the $C^{1+\alpha}(\mathbb{R}^n)$ should be optimal.
  \end{remark}
Like the classical case, we can get the $C^{2+\alpha}$-estimate for the $x$-axis if
$f(\cdot,y)\in C_b^\beta(\mathbb{R})$ for
any $y\in\mathbb{R}$, but we can not get the same estimate for $y$-axis. If we
want to get the $C^{2+\alpha}$-estimate for the $y$-axis, we need more regularity
about the second variable. In other words, we have the following results.
\begin{col}\lbl{c0.1}  Assume that $f(\cdot,y)\in C_b^\beta(\mathbb{R})$ for
any $y\in\mathbb{R}$ and
$\nabla_yf(x,\cdot)\in C_b^\beta(\mathbb{R})$ for any $x\in\mathbb{R}$. Then
  \bess
[\nabla^2u]_\beta\leq C\lambda^{-\delta}\|f\|_{C^\beta(\mathbb{R};C^{1+\beta}(\mathbb{R}))},
  \eess
where $0<\delta<(\beta\wedge(1-\beta))$. That is to say, $\nabla^2 u\in C_b^\beta(\mathbb{R}^2)$.
  \end{col}

{\bf Proof of Theorem \ref{t1.1}.} We use Picard's iteration to solve (\ref{0.1}). Let $u_0=0$ and define for $n\in \mathbb{N}$,
   \bes
u_n:=(\lambda-\mathcal{L})^{-1}(f-b\cdot\nabla u_{n-1}).
   \lbl{0.3}\ees
It follows from Lemma \ref{l2.1} that
   \bes
\lambda^{\delta}\|u_n\|_{C_b^{1+\beta}(\mathbb{R}^2)}&\leq&
\|f-b\cdot\nabla u_{n-1}\|_{ C_b^{\beta}(\mathbb{R}^2)}\nm\\
&\leq& \|f\|_{ C_b^{\beta}(\mathbb{R}^2)}
+\|b\|_{L^\infty(\mathbb{R}^2)}\cdot\|\nabla u_{n-1}\|_{ C_b^{\beta}(\mathbb{R}^2)}\nm\\
&&
+\|b\|_{ C_b^{\beta}(\mathbb{R}^2)}\cdot\|\nabla u_{n-1}\|_{L^\infty(\mathbb{R}^2)},
  \lbl{0.4}\ees
and
   \bes
\lambda^{\delta}\|u_n-u_m\|_{C_b^{1+\beta}(\mathbb{R}^2)}
&\leq& \|b\|_{L^\infty(\mathbb{R}^2)}\cdot\|\nabla u_{n-1}-\nabla u_{m-1}\|_{ C_b^{\beta}(\mathbb{R}^2)}\nm\\
&&
+\|b\|_{ C_b^{\beta}(\mathbb{R}^2)}\cdot\|\nabla u_{n-1}-\nabla u_{m-1}\|_{L^\infty(\mathbb{R}^2)}.
  \lbl{0.5}\ees
Choosing $\lambda_0$ be large enough so that $C\lambda^{-\delta}\|b\|_{C_b^\beta({\mathbb{R}^2})}<1/4$ for
all $\lambda\geq\lambda_0$, we get
   \bess
\|u_n\|_{C_b^{1+\beta}(\mathbb{R}^2)}\leq\lambda^{-\delta}(\|f\|_{ C_b^{\beta}(\mathbb{R}^2)}
+2\|b\|_{ C_b^{\beta}(\mathbb{R}^2)})
+\frac{1}{2}\| u_{n-1}\|_{ C_b^{1+\beta}(\mathbb{R}^2)}
   \eess
and for all $n\geq m$,
   \bess
\|u_n-u_m\|_{C_b^{1+\beta}(\mathbb{R}^2)}\leq \frac{1}{2}\|u_{n-1}-u_{m-1}\|_{C_b^{1+\beta}(\mathbb{R}^2)}.
  \eess
Substituting them into (\ref{0.4}) and (\ref{0.5}), we obtain
  \bess
\lambda^{-\delta}\|u_n\|_{C_b^{1+\beta}(\mathbb{R}^2)}\leq C\|f\|_{ C_b^{\beta}(\mathbb{R}^2)}
   \eess
and for all $n\geq m$,
   \bess
\lambda^{-\delta}\|u_n-u_m\|_{C_b^{1+\beta}(\mathbb{R}^2)}\leq \frac{C}{2^m}.
   \eess
Hence there is a $u\in C_b^{1+\beta}(\mathbb{R}^2)$ such that (\ref{0.2}) holds and
  \bess
\lambda^{-\delta}\|u-u_m\|_{C_b^{1+\beta}(\mathbb{R}^2)}\leq \frac{C}{2^m},
   \eess
and $u$ solves equation (\ref{0.1}) by taking limits for (\ref{0.3}). The
second result can be obtained similarly. The proof is complete. $\Box$

{\bf Proof of Theorem \ref{t1.2}.} For simplicity, we only
consider a special case, that is, $f(x,y)=f_1(x)f_2(y)$. Assume that $f_1\in L^p(\mathbb{R})$
and $f_2\in L^q(\mathbb{R})$.
Denote
   \bess
K_1(t,y)=\int_{\mathbb{R}}K(t,x-u,y)f_1(u)du.
   \eess
By using the above inequality, Minkowski's inequality and the properties of
the heat kernel, we have
 \bess
\|u\|_{L^r(\mathbb{R}^2)}&=&\|\int_0^\infty e^{-\lambda t}K(t,\cdot)\ast f(x,y)dt\|_{L^r(\mathbb{R}^2)}\\
&\leq&\int_0^\infty e^{-\lambda t}\left(\int_{\mathbb{R}^2}|\int_{\mathbb{R}}K_1(t,x,y-v) f_2(v)dv|^rdxdy\right)^{\frac{1}{r}}dt\\
&=&\int_0^\infty e^{-\lambda t}\left(\int_{\mathbb{R}}\|K_1(t,x,\cdot)\ast f_2\|^r_{L^r(\mathbb{R})}dx\right)^{\frac{1}{r}}dt\\
&\leq&\|f_2\|_{L^q(\mathbb{R})}\int_0^\infty e^{-\lambda t}\left(\int_{\mathbb{R}}\|K_1(t,x,\cdot)\|^r_{L^m(\mathbb{R})}dx\right)^{\frac{1}{r}}dt\\
&=&\|f_2\|_{L^q(\mathbb{R})}\int_0^\infty e^{-\lambda t}\left(\int_{\mathbb{R}}\Big|\int_\mathbb{R}|K_1(t,x,y)|^mdy\Big|^{\frac{r}{m}}dx\right)^{\frac{1}{r}}dt\\
&\leq&\|f_2\|_{L^q(\mathbb{R})}\int_0^\infty e^{-\lambda t}\left(\int_{\mathbb{R}}\Big|\int_\mathbb{R}|\int_{\mathbb{R}}K(t,x-u,y)f_1(u)du|^rdx\Big|^{\frac{m}{r}}dy\right)^{\frac{1}{m}}dt\\
&=&\|f_2\|_{L^q(\mathbb{R})}\int_0^\infty e^{-\lambda t}\left(\int_{\mathbb{R}}\Big|\int_\mathbb{R}\|K(t,\cdot,y)\ast f_1\|^m_{L^r(\mathbb{R})}dy\right)^{\frac{1}{m}}dt\\
&\leq&\|f_1\|_{L^p(\mathbb{R})}\|f_2\|_{L^q(\mathbb{R})}\int_0^\infty e^{-\lambda t}\|K(t,\cdot,\cdot)\|_{L^{n}(\mathbb{R},L^{m}(\mathbb{R}))}dt\\
&\leq&C_1\int_0^\infty e^{-\lambda t}t^{-\frac{3}{2}+\frac{1}{2n}+\frac{1}{m}}dt\\
&\leq&C\lambda^{-1-\frac{3}{2r}+\frac{1}{2p}+\frac{1}{q}},
  \eess
where
   \bes
&&1+\frac{1}{r}=\frac{1}{n}+\frac{1}{p}, \ \ \ \ 1+\frac{1}{r}=\frac{1}{m}+\frac{1}{q},\lbl{3.1}\\
&&C_1=\|f_1\|_{L^p(\mathbb{R})}\|f_2\|_{L^q(\mathbb{R})}
\left(\int_{\mathbb{R}}\left(\int_{\mathbb{R}}e^{-p'x^2-p'(y-x^2)^2}dx\right)^{\frac{q'}{p'}}dy\right)^{\frac{1}{q'}}.\nm
  \ees
We remark that
   \bess
&&\left(\int_\mathbb{R}|\nabla_xK(t,x,y)|^rdx\right)^{\frac{1}{r}}\leq Ct^{-\frac{1}{2}+\frac{1}{2r}},\\
&&\left(\int_\mathbb{R}|\nabla_yK(t,x,y)|^rdy\right)^{\frac{1}{r}}\leq Ct^{-1+\frac{1}{r}}.
  \eess
Similarly, we obtain
 \bess
\|\nabla_xu\|_{L^r(\mathbb{R}^2)}&\leq&C\lambda^{-\frac{1}{2}-\frac{3}{2r}+\frac{1}{2p}+\frac{1}{q}}.
   \eess
Furthermore, we can get
 \bess
\|\nabla^2_xu\|_{L^r(\mathbb{R}^2)}\leq C\lambda^{-\frac{3}{2r}+\frac{1}{2p}+\frac{1}{q}}.
   \eess
Thus if we take $p=q=r$, then we have $\|\nabla^2_xu\|_{L^r(\mathbb{R}^2)}\leq C$.
However, if we deal with the second variable, it is difficult to get the decay estimate. More precisely,
we have for $\frac{3}{2r}<\frac{1}{2p}+\frac{1}{q}$,
 \bess
\|\nabla_yu\|_{L^r(\mathbb{R}^2)}&\leq&C\lambda^{-\frac{3}{2r}+\frac{1}{2p}+\frac{1}{q}}\to\infty, \ \ {\rm as}\ \ \lambda\to\infty.
   \eess
Hence we must add more regularity on the second variable. Meanwhile, we recall that if $h\in W^{s,p}(\mathbb{R}^n)$ with $0<s<1$, then
  \bess
\|h\|_{W^{s,p}(\mathbb{R}^n)}=\left(\|h\|^p_{L^p(\mathbb{R}^n)}+
\int_{\mathbb{R}^n}\int_{\mathbb{R}^n}\frac{|h(x)-h(y)|^p}{|x-y|^{n+sp}}dxdy\right)^{\frac{1}{p}}.
   \eess
If we assume that $f_2\in W^{s,q}(\mathbb{R})$ and let
   \bess
K_2(t,y)=\int_{\mathbb{R}}\nabla_yK(t,x,y-v)f_2(v)dv.
   \eess
 we get
 \bess
\|\nabla_yu\|_{L^r(\mathbb{R}^2)}&=&\|\int_0^\infty e^{-\lambda t}\nabla_yK(t,\cdot)\ast f(x,y)dt\|_{L^r(\mathbb{R}^2)}\\
&\leq&\int_0^\infty e^{-\lambda t}\left(\int_{\mathbb{R}^2}|\int_{\mathbb{R}}K_2(t,x-u,y) f_1(u)du|^rdxdy\right)^{\frac{1}{r}}dt\\
&\leq&\|f_1\|_{L^p(\mathbb{R})}\int_0^\infty e^{-\lambda t}\left(\int_{\mathbb{R}}\|K_2(t,\cdot,y)\|^n_{L^r(\mathbb{R})}dy\right)^{\frac{1}{n}}dt\\
&=&\|f_1\|_{L^p(\mathbb{R})}\int_0^\infty e^{-\lambda t}\\
&&\times\left(\int_\mathbb{R}\left(\Big|\int_\mathbb{R}\nabla_yK(t,x,v)|v|^{s+\frac{1}{q}}\frac{f_2(y-v)-f_2(y)}{|v|^{s+\frac{1}{q}}}dv\Big|^rdx
\right)^{\frac{n}{r}}dy\right)^{\frac{1}{n}}dt\\
&\leq&\|f_1\|_{L^p(\mathbb{R})}\|f_2\|_{W^{s,q}(\mathbb{R})}\int_0^\infty e^{-\lambda t}\left(\int_\mathbb{R}\left(|\nabla_yK(t,x,y)|^m|y|^{sm+\frac{m}{q}}dy\right)^{\frac{n}{m}}dx\right)^{\frac{1}{n}}dt\\
&\leq&C_2\int_0^\infty e^{-\lambda t}t^{-\frac{5}{2}+s+\frac{1}{q}+\frac{1}{2n}+\frac{1}{m}}dt\\
&\leq&C\lambda^{-s-\frac{3}{2r}+\frac{1}{2p}},
  \eess
where $m,n$ satisfy (\ref{3.1}) and
   \bess
C_2=2\|f_1\|_{L^p(\mathbb{R})}\|f_2\|_{W^{s,q}(\mathbb{R})}
\left(\int_{\mathbb{R}}\left(\int_{\mathbb{R}}|y-x^2|^p|y|^{sp+\frac{p}{q}}e^{-p'x^2-p'(y-x^2)^2}dx\right)^{\frac{q'}{p'}}dy\right)^{\frac{1}{q'}}.
  \eess
Moreover,
under the condition that $f_2\in W^{s,p}(\mathbb{R})$, we can similarly get
 \bess
\|\nabla^2_yu\|_{L^r(\mathbb{R}^2)}&\leq&C\lambda^{-s-\frac{3}{2r}+\frac{1}{2p}+1}.
   \eess
Hence it is easy to see that we can take suitable $s\in(0,1),p>1$ and $r>1$ such that
$-s-\frac{3}{2r}+\frac{1}{2p}+1\leq0$. That is to say, we have
 \bess
\|\nabla^2_yu\|_{L^r(\mathbb{R}^2)}&\leq&C
   \eess
under the condition that $f_2\in W^{s,p}(\mathbb{R})$. The proof  is complete. $\Box$

\begin{remark}\lbl{r2.2} It is well known that if $f\in L^p(\mathbb{R}^n)$, then the solution $u$ of the following
equation
   \bess
u_t-\Delta u=f,\ \ u_0=0,
   \eess
belongs to $W^{2,p}(\mathbb{R}^n)$, which is the $L^p$-theory. Noting that the heat kernel of above
equation is Gauss heat kernel, and similar to Remark \ref{r2.1}, it is easy to find the 
difference from the classical Laplacian operator. Due to the singularity of the variable $y$, we 
must give two different assumptions. Comparing the classical $L^p$-theory, 
in our case the regularity of Theorem \ref{t1.2} should be optimal.
  \end{remark}

\medskip

\noindent {\bf Acknowledgment} The first author was supported in part
by NSFC of China grants 11571176 and 11771123.

 \end{document}